\magnification=1200  \voffset=1cm  \hoffset=1cm
\hsize=11.25cm
\vsize=18cm
\parskip 0pt
\parindent=12pt



\catcode'32=9

\font\tenpc=cmcsc10
\font\eightpc=cmcsc8
\font\eightrm=cmr8
\font\eighti=cmmi8
\font\eightsy=cmsy8
\font\eightbf=cmbx8
\font\eighttt=cmtt8
\font\eightit=cmti8
\font\eightsl=cmsl8
\font\sixrm=cmr6
\font\sixi=cmmi6
\font\sixsy=cmsy6
\font\sixbf=cmbx6

\skewchar\eighti='177 \skewchar\sixi='177
\skewchar\eightsy='60 \skewchar\sixsy='60

\catcode`@=11

\def\tenpoint{%
  \textfont0=\tenrm \scriptfont0=\sevenrm \scriptscriptfont0=\fiverm
  \def\rm{\fam\z@\tenrm}%
  \textfont1=\teni \scriptfont1=\seveni \scriptscriptfont1=\fivei
  \def\oldstyle{\fam\@ne\teni}%
  \textfont2=\tensy \scriptfont2=\sevensy \scriptscriptfont2=\fivesy
  \textfont\itfam=\tenit
  \def\it{\fam\itfam\tenit}%
  \textfont\slfam=\tensl
  \def\sl{\fam\slfam\tensl}%
  \textfont\bffam=\tenbf \scriptfont\bffam=\sevenbf
  \scriptscriptfont\bffam=\fivebf
  \def\bf{\fam\bffam\tenbf}%
  \textfont\ttfam=\tentt
  \def\tt{\fam\ttfam\tentt}%
  \abovedisplayskip=12pt plus 3pt minus 9pt
  \abovedisplayshortskip=0pt plus 3pt
  \belowdisplayskip=12pt plus 3pt minus 9pt
  \belowdisplayshortskip=7pt plus 3pt minus 4pt
  \smallskipamount=3pt plus 1pt minus 1pt
  \medskipamount=6pt plus 2pt minus 2pt
  \bigskipamount=12pt plus 4pt minus 4pt
  \normalbaselineskip=12pt
  \setbox\strutbox=\hbox{\vrule height8.5pt depth3.5pt width0pt}%
  \let\bigf@ntpc=\tenrm \let\smallf@ntpc=\sevenrm
  \let\petcap=\tenpc
  \normalbaselines\rm}

\def\eightpoint{%
  \textfont0=\eightrm \scriptfont0=\sixrm \scriptscriptfont0=\fiverm
  \def\rm{\fam\z@\eightrm}%
  \textfont1=\eighti \scriptfont1=\sixi \scriptscriptfont1=\fivei
  \def\oldstyle{\fam\@ne\eighti}%
  \textfont2=\eightsy \scriptfont2=\sixsy \scriptscriptfont2=\fivesy
  \textfont\itfam=\eightit
  \def\it{\fam\itfam\eightit}%
  \textfont\slfam=\eightsl
  \def\sl{\fam\slfam\eightsl}%
  \textfont\bffam=\eightbf \scriptfont\bffam=\sixbf
  \scriptscriptfont\bffam=\fivebf
  \def\bf{\fam\bffam\eightbf}%
  \textfont\ttfam=\eighttt
  \def\tt{\fam\ttfam\eighttt}%
  \abovedisplayskip=9pt plus 2pt minus 6pt
  \abovedisplayshortskip=0pt plus 2pt
  \belowdisplayskip=9pt plus 2pt minus 6pt
  \belowdisplayshortskip=5pt plus 2pt minus 3pt
  \smallskipamount=2pt plus 1pt minus 1pt
  \medskipamount=4pt plus 2pt minus 1pt
  \bigskipamount=9pt plus 3pt minus 3pt
  \normalbaselineskip=9pt
  \setbox\strutbox=\hbox{\vrule height7pt depth2pt width0pt}%
  \let\bigf@ntpc=\eightrm \let\smallf@ntpc=\sixrm
  \let\petcap=\eightpc
  \normalbaselines\rm}
\catcode`@=12

\tenpoint



\font\tengoth=eufm10


\catcode`\@=11
\def\pc#1#2|{{\bigf@ntpc #1\penalty \@MM\hskip\z@skip\smallf@ntpc%
	\uppercase{#2}}}
\catcode`\@=12

\def\pointir{\discretionary{.}{}{.\kern.35em---\kern.7em}\nobreak
   \hskip 0em plus .3em minus .4em }

\def\qed{\quad\raise -2pt\hbox{\vrule\vbox to 10pt{\hrule width 4pt
   \vfill\hrule}\vrule}}

\def\cqfd{\penalty 500 \hbox{\qed}\par\bigskip}

\def\rem#1|{\par\medskip{{\it #1}\pointir}}

\def\vspace[#1]{\noalign{\vskip#1}}

\def\abstract#1{\vbox{\eightpoint\narrower\narrower
\pc ABSTRACT|\pointir #1}}


\def\enslettre#1{\font\zzzz=msbm10 \hbox{\zzzz #1}}
\def\smallenslettre#1{\font\zzzz=msbm7 \hbox{\zzzz #1}}

\def\setZ{{\enslettre Z}}
\def\setQ{{\enslettre Q}}

\def\setC{{\enslettre C}}
\def\setF{{\enslettre F}}
\def\setQsmall{{\smallenslettre Q}}
\def\setCsmall{{\smallenslettre C}}

\long\def\maskbegin#1\maskend{}

\def\section#1{\goodbreak\par\vskip .5cm\centerline{\bf #1}
   \par\nobreak\vskip 5pt }

\long\def\th#1|#2\endth{\par\medbreak
   {\petcap #1\pointir}{\it #2}\par\medbreak}

\def\article#1|#2|#3|#4|#5|#6|#7|
    {{\leftskip=7mm\noindent
     \hangindent=2mm\hangafter=1
     \llap{{\tt [#1]}\hskip.35em}{\petcap#2}\pointir
     #3, {\sl #4}, {\bf #5} ({\oldstyle #6}),
     pp.\nobreak\ #7.\par}}
\def\livre#1|#2|#3|#4|
    {{\leftskip=7mm\noindent
    \hangindent=2mm\hangafter=1
    \llap{{\tt [#1]}\hskip.35em}{\petcap#2}\pointir
    {\sl #3}, #4.\par}}
\def\divers#1|#2|#3|
    {{\leftskip=7mm\noindent
    \hangindent=2mm\hangafter=1
     \llap{{\tt [#1]}\hskip.35em}{\petcap#2}\pointir
     #3.\par}}


\font\KFracFont=cmr12 at 20pt
\def\KK{\mathop{\lower 4pt\hbox{\KFracFont K}}\limits}

\overfullrule=0pt


\def\rem{\noindent {\bf{Remark : }}}

\def\house#1{\setbox1=\hbox{$\,#1\,$}%
\dimen1=\ht1 \advance\dimen1 by 2pt \dimen2=\dp1 \advance\dimen2 by 2pt
\setbox1=\hbox{\vrule height\dimen1 depth\dimen2\box1\vrule}%
\setbox1=\vbox{\hrule\box1}%
\advance\dimen1 by .4pt \ht1=\dimen1
\advance\dimen2 by .4pt \dp1=\dimen2 \box1\relax}

  \def\deg{{\rm deg}}
   
  \def\eps{{\varepsilon}}
       
   \def\th{{\theta}}

\def\build#1_#2^#3{\mathrel{\mathop{\kern 0pt#1}\limits_{#2}^{#3}}}

\def\date {le\ {\the\day}\ \ifcase\month\or
janvier\or fevrier\or mars\or avril\or mai\or juin\or juillet\or
ao\^ut\or septembre\or octobre\or novembre\or
d\'ecembre\fi\ {\oldstyle\the\year}}

\font\fivegoth=eufm5 \font\sevengoth=eufm7 \font\tengoth=eufm10

\newfam\gothfam \scriptscriptfont\gothfam=\fivegoth
\textfont\gothfam=\tengoth \scriptfont\gothfam=\sevengoth

\def\cqfd{\unskip\kern 6pt\penalty 500 \raise 0pt\hbox{\vrule\vbox
to6pt{\hrule width 6pt \vfill\hrule}\vrule}\par}

\def\smallsquare{\vbox{\hrule\hbox{\vrule height 1 ex\kern 1 ex\vrule}\hrule}}


\vskip 5mm

\def\qed{\quad\raise -2pt\hbox{\vrule\vbox to 10pt{\hrule width 4pt
\vfill\hrule}\vrule}}

\def\cqfd{\qed}

\def\pointir{\discretionary{.}{}{.\kern.35em---\kern.7em}\nobreak
\hskip 0em plus .3em minus .4em }



\bigskip

\centerline{\bf Hankel determinants, Pad\'e approximations,}
\centerline{\bf and irrationality exponents}


\bigskip
\centerline{\sl Yann BUGEAUD, Guo-Niu HAN, Zhi-Ying WEN, Jia-Yan YAO}
\footnote{}{\eightpoint
}
\bigskip
\bigskip
{\narrower\narrower
\eightpoint
\noindent
{\bf Abstract}.\quad
The irrationality exponent of an irrational number $\xi$, 
which measures the approximation rate
of $\xi$ by rationals,  is in general extremely    
difficult to compute explicitly, unless we know the continued fraction expansion of $\xi$. 
Results obtained so far 
are rather fragmentary and often treated case by case. 
In this work, we shall unify all the known
results on the subject by showing that the irrationality exponents of large classes of
automatic numbers and Mahler numbers (which are transcendental) are exactly equal to $2$.
Our classes contain the Thue--Morse--Mahler numbers,
the sum of the reciprocals of the Fermat numbers, the regular paperfolding numbers,
which have been previously considered respectively by Bugeaud, Coons, and Guo, Wu and Wen,
but also new classes such as the Stern numbers and so on.
Among other ingredients, our proofs use results on Hankel determinants obtained recently by Han.
}

\section{1. Introduction} 

Let $\xi$ be an irrational real number. The irrationality exponent
$\mu(\xi)$ of $\xi$ is the supremum
of the real numbers $\mu$ such that the inequality
$$
\biggl| \xi - {p \over q} \biggr| < {1 \over q^{\mu}}
$$
has infinitely many solutions in rational numbers $p/q$. Hence, we have
$$
\mu(\xi)=1-\liminf\limits_{q\rightarrow \infty}{\log\|q\xi\|\over\log q},\leqno{(1.1)}
$$
where $\|x\|$ denotes the distance between the real number $x$ and its nearest integer.
An easy covering argument shows that
$\mu (\xi)$ is at most  
equal to $2$ for almost
all real numbers $\xi$ (with respect
to the Lebesgue measure). It follows from the theory of continued fractions
that the irrationality exponent of an irrational real number is always greater than or equal
to $2$. More precisely, let $[a_0; a_1,a_2,\ldots]$ denote the continued fraction expansion
of an irrational real number $\xi$ and $(p_n / q_n)_{n \ge 1}$ denote the
sequence of its convergents
(for more about continued fractions, see for example [La95]). Then, we have
$$
\mu(\xi)=2+\limsup_{n\rightarrow\infty}{\log a_{n+1}\over \log q_n}.\leqno{(1.2)}
$$

Furthermore, Roth's theorem [Ro55]
asserts that the irrationality exponent of every
algebraic irrational number is equal to $2$.
However, it is in general a very difficult problem to determine the
irrationality exponent of a given transcendental real number $\xi$.
Apart from some numbers involving the exponential
function or the Bessel function (see the end of Section 1
in [Ad10]) and apart from
more or less {\it ad hoc} constructions (see below),
it seems to us that
the only known method
to determine the irrationality
exponent of (certain) transcendental numbers is the method developed in [Bu11].
Up to now, this method has been applied to a handful of irrational
numbers [Bu11, Co13, GWW14, WW14]. The main purpose of the present work
is to considerably extend these results and to exhibit infinite families of
transcendental numbers with irrationality exponent equal to $2$.

Let us now focus on a special class of real numbers.

A real number $\xi$ is automatic if there exist two integers $k,b\geq 2$
such that the $b$-ary expansion of $\xi$ is $k$-automatic.
This means that, if we write
$\xi=\sum\limits_{n\geq 0}{a(n)\over b^n}$ with $a(n)\in \setZ\ (n\geq 0)$
and $0\leq a(n)<b$ for $n\geq 1$, then the set of subsequences
$$
\Big\{\big(a(k^rn+s)\big)_{n\geq 0}\ |\ r\geq 0,\ 0\leq s<k^r\Big\}
$$
is finite (For more on automatic sequences,
see for example Allouche [Al87] and also the book of Allouche and Shallit [AS03]).
For example, the case of Kmo\v{s}ek-Shallit
numbers $f_{KS}\big({1\over b}\big)=\sum\limits_{n\geq 0}{1\over b^{2^n}}$
(studied independently by Kmo\v{s}ek [Km79] and Shallit [Sh79] in 1979
to give ``natural'' examples of real numbers with bounded partial quotients)
corresponds
to the characteristic function of the set $\{2^n\ | \ n\geq 0\}$,
which is $2$-automatic but not ultimately periodic. These numbers are  
transcendental (see [Ke16], [Ma29], and also [LVdP77]).
It was long conjectured and finally has been proved by Adamczewski and Bugeaud
in [AB07] that an automatic number is either rational or transcendental;
see [BBC15] and [Ph15] for two recent alternative proofs.  
We have thus a large family of
``simple'' transcendental numbers, and one can then ask what are their irrationality exponents.

In 2006, Adamczewski and Cassaigne showed in [AC06]
that an automatic number cannot be a Liouville number (recall that, by definition,
a Liouville number is a real number whose irrationality exponent is infinite).
Subsequently, Adamczewski and Rivoal [AR09]  obtained in 2009
upper bounds for the irrationality exponents
of some famous automatic numbers constructed from
the Thue--Morse, Rudin--Shapiro, paperfolding and Baum--Sweet sequences.
In 2008, Bugeaud [Bu08]  constructed explicitly elements of the classical middle third Cantor
set with any prescribed irrationality exponent
(an analog for the function field case has been obtained very recently by Pedersen [Pe14]),
and proved that there exist automatic real numbers
with any prescribed rational irrationality exponent.
But what is the exact value of the irrationality exponent of a given automatic
irrational number (for example, the famous Thue--Morse--Mahler numbers)?
This question was addressed in [BKS11], and the results obtained on this subject
are rather fragmentary even until now,
often treated case by case, and can be summarized as follows.

The history begun in 2011 with the paper [Bu11], in which Bugeaud
developed a method to show that the irrationality exponents of the
Thue--Morse--Mahler numbers are equal to $2$.
Recall that the famous Thue--Morse sequence $(t_n)_{n\geq 0}$ on $\{0,1\}$
is defined recursively by $t_0=0$, $t_{2n}=t_n$ and $t_{2n+1}=1-t_n$ for all
integers $n\geq 0$, and that the Thue--Morse--Mahler numbers take the form
$$
f_{TMM} \Bigl({1\over b} \Bigr)=\sum\limits_{n\geq 0}{t_n\over b^n},
$$
where $b \ge 2$ is an integer.
Recall also that the Thue--Morse
sequence is $2$~-automatic but not ultimately periodic, and Mahler [Ma29] already showed in 1929
that $f_{TMM}(1/2)$ is transcendental (see also Dekking [De77] for another proof).

In 2013, Coons considered in [Co13] the following two power series
$$
{\cal F}(z)=\sum\limits_{n\geq 0}{{z^{2^n}}\over{1+z^{2^n}}},\quad
{\cal G}(z)=\sum\limits_{n\geq 0}{{z^{2^n}}\over{1-z^{2^n}}},\leqno (1.3)
$$
and showed that for all integers $b\geq 2$,
we have $\mu({\cal F}(1/b))=\mu({\cal G}(1/b))=2$.
Note here that the special value
${\cal F}(1/2)$ is the sum of the reciprocals of the Fermat numbers $F_n:=2^{2^n}+1$,
and the sequence of coefficients of ${\cal G}(z)$
is usually called the {\it Gros sequence} [Gr72, HKMP13].

In 2014, Guo, Wu and Wen considered in [GWW14]
the regular paperfolding numbers defined by
$$
f_{RPF}\Big({1\over b}\Big):=\sum\limits_{n\geq 0}{{u_n}\over{b^n}},
$$
where $b\geq 2$ is an integer, and $(u_n)_{n\geq 0}$
is the regular paperfolding sequence on $\{0,1\}$ defined recursively by
$u_{4n}=1$, $u_{4n+2}=0$, and $u_{2n+1}=u_n$, for~all integers $n\geq 0$.
They proved that the irrationality exponents of these numbers are all equal to~$2$.
For more on the regular paperfolding sequence, see for example [Al87] and [AS03].

Very recently, Wen and Wu [WW14] studied the Cantor real numbers
$$
f_{C}\Big({1\over b}\Big):=\sum\limits_{n\geq 0}{{v_n}\over{b^n}},
$$
where $b\geq 2$ is an integer, and $(v_n)_{n\geq 0}$
is the Cantor sequence on $\{0,1\}$ such that for all integers $n\geq 0$,
we have $v_n=1$ if and only if the ternary expansion of $n$ does not contain the digit $1$.
They showed that
the irrationality exponents of these numbers are also equal to $2$.
We point out that the Cantor sequence is $3$-automatic (see for example [AS03]) and that
its generating function $f_C$ satisfies $f_C (z) = (1 + z^2) f_C (z^3)$.

In the present work, we shall unify all the above results together
and compute the irrationality exponent of some
new families of transcendental numbers.
We do not restrict our attention to automatic numbers and take a more general point of view.

Mahler's method [Ma29, Ma30a, Ma30b] is a method in transcendence theory
whereby one uses a function $F(z) \in \setQ[[z]]$ that satisfies a functional equation
of the following form
$$
\sum_{i=0}^n P_i(z) F(z^{d^i}) = 0,     \leqno (1.4)
$$
for some integers $n \geq 1$ and $d\geq 2$, and polynomials $P_0(z),
\ldots, P_n (z)$ in~$\setZ[x]$
with $P_0 (z) P_n(z)\neq 0$, to give results about the nature of the numbers
$F(1/b)$ with $b\geq 2$ an integer such that $1/b$ is less
than the radius of convergence of $F(z)$.
We refer to such numbers $F(1/b)$ as {\it Mahler numbers}.
It is well known that automatic numbers
are special cases of Mahler numbers (see [Be94, Theorem 1]).
The following theorem, established in [BBC15], extends the main result
of [AC06], quoted above.

\proclaim Theorem 1.1.
A Mahler number cannot be a Liouville number.

By means of a suitable adaptation of the so-called Mahler's method, 
it is proved in [BBC15] that an irrational Mahler number is
transcendental
when $P_0 (z)$ in (1.4) is a nonzero integer. However, the general case
remains an open problem. 
Note that Corvaja and Zannier [CZ02] explained how,  
beside Mahler's method, the Schmidt Subspace Theorem can be used to prove,  
under quite general assumptions, the transcendence of values of power series 
with integer coefficients at non-zero algebraic points.     

We formulate the following open question.

\proclaim Problem 1.2.
To determine the set of irrationality exponents of
irrational Mahler numbers.

Actually we will consider power series $F(z)$ satisfying a functional equation
of the special form
$$
P_{-1} (z) + P_0 (z) F(z) + P_1(z) F(z^{d}) = 0,     \leqno (1.5)
$$
for some integer $d\geq 2$ and polynomials $P_{-1} (z), P_0(z), P_1 (z) \in \setZ[x]$, 
with $P_0 (z)$ and $P_1(z)$ being non-zero. 
Observe that, by combining (1.5) with the equation obtained by
substituting $z$ with $z^d$ in (1.5), we see that $F(z)$ satisfies an equation of the type (1.4).
We also point out here that by a general result of Zannier [Za98, p.~18],  the function $F(z)$
is either rational or transcendental over $\setQ(z)$.

The present work is organized as follows. In Section 2, we highlight several of our results.
Then, in Section 3, we recall some basic
notation and results about Pad\'e approximation, which is the starting point
of our study. In Section 4, we compute with Hankel determinants
the irrationality exponent of certain transcendental numbers, which
are values at the inverse of integers $\geq 2$ of power
series satisfying a functional equation of type~(1.5).
Since it is extremely difficult to compute explicitly the Hankel
determinants of a given sequence, we collect, in Section 5, some results
about Hankel continued fractions obtained very recently by Han [H15a, H15b], and apply them
in Section 6 to obtain directly (this means, without condition on Hankel determinants)   
the irrationality exponent of special values of some power series
satisfying a special type of functional equation.
Our results cover all the known results on irrationality exponent
listed above, and in the final Section 7, we shall give several
new applications to obtain the irrationality exponent of 
new families of transcendental numbers.


\section{2. Results} 

Let $d\geq 2$ be an integer, and $(c_m)_{m\geq 0}$
be an integer sequence such that
$f(z)=\sum\limits_{m=0}^{+\infty}c_m z^m$ converges inside the unit disk.
Suppose that there exist integer polynomials $A(z)$, $B(z)$, $C(z)$, and $D(z)$ such that
$$
f(z)={A(z)\over B(z)} + {C(z)\over D(z)} f(z^d).        \leqno{(2.1)}
$$
Under various assumptions on these polynomials, we are able to show that, for every
integer $b \ge 2$, the irrationality exponent of $f(1/b)$ is equal to $2$.

One of our tools is a careful study of the sequence $(H_n (f))_{n \ge 0}$
of the Hankel determinants of $f$, defined by $H_0 (f) = 1$ and
$$
H_n (f) := \left|
\matrix{ c_0 & c_{1} & \ldots & c_{n-1} \cr
c_{1} & c_{2} & \ldots & c_{n} \cr
\ \vdots \hfill & \ \vdots \hfill & \ddots &
\ \vdots \hfill \cr
c_{n-1} & c_{n} & \ldots & c_{2n-2} \cr} \right| , \quad \hbox{for all integers $n \ge 1$}.
$$

We state below a consequence of our Theorem 4.1, which highlights a relationship
between the irrationality exponent of $f(1/b)$ and the sequence $(H_n (f))_{n \ge 0}$,
and correct, improve and generalize the main result
recently obtained by Guo, Wu, and Wen [GWW14].

\proclaim Theorem 2.1.
Let $d\geq 2$ be an integer, and $(c_j)_{j\geq 0}$ be an integer sequence such that
$f(z)=\sum\limits_{j=0}^{+\infty}c_j z^j$ converges inside the unit disk.
Suppose that there exist integer polynomials $A(z)$, $B(z)$ and $C(z)$ such that
$$
f(z)={A(z)\over B(z)} + C(z)f(z^d).
$$
Let $b \ge 2$ be an integer such that
$C({1\over b^{d^m}})\neq 0$
for all integers $m\geq 0$.
If there exists an increasing sequence of positive integers
$(n_i)_{i\geq 0}$  such that $H_{n_i}(f)\not=0$ for all integers $i\geq 0$ and
$\lim\limits_{i\rightarrow \infty} {n_{i+1}\over n_i} =1$,
then $f(1/b)$ is transcendental and its irrationality exponent is equal to $2$.
\vskip 5pt

\noindent {\bf Remark.}
In [GWW14] the authors need to assume the existence of an infinite sequence
$(n_i)_{i \ge 1}$ satisfying
$\liminf\limits_{i\rightarrow \infty}{{n_{i+1}}\over {n_i}}=1$
and such that $H_{n_i} (f) H_{n_i + 1} (f)$ is nonzero for all integers $i \ge 1$.
However, in their proof, they
make use of the stronger assumption that this limit inferior is actually
a limit, and also use implicitly the fact that $C({1\over b^{d^m}})\neq 0$
for all integers $m\geq 0$.
Thus, our general result Theorem 4.1 considerably extends (and corrects)
Theorem~1 of [GWW14].
\vskip 5pt

Theorem 2.1 will be proved in Section 4.
\vskip 5pt

However, the computation of the sequence $(H_n (f))_{n \ge 0}$ is not an easy task,
and even to get information on its vanishing terms is difficult. Very recently,
Han [H15a, H15b] has developed a new and fruitful method. As a result, we obtain in particular
the following theorem.

\proclaim Theorem 2.2.
Let $f(z)\in \setZ[[z]]$ be the power series defined by
$$
f(z)= \prod_{n\geq 0}\Bigl(1+uz^{2^n}+2z^{2^{n+1}} {C(z^{2^n})\over D(z^{2^n})}\Bigr), \leqno{(2.2)}
$$
where $u\in\setZ$, and $C(z), D(z)\in \setZ[z]$ with
$D(0)= 1$. Let $b\geq 2$ be an integer such that
$D\big({1\over b^{2^m}}\big)f\big({1\over b^{2^m}}\big)\neq 0$ for all integers $m\geq 0$.
If $f(z)\, ({\rm mod}\, 4)$ is not a rational function,
then $f(1/b)$ is transcendental and its irrationality exponent is equal to $2$.

\noindent {\bf Remark.} Taking $C(z)=0$, $D(z)=1$, and $u=-1$ in Theorem 2.2,
we recover the result of [Bu11] about Thue--Morse--Mahler numbers
which states that $\mu(f_{TMM}(1/b))=2$, for all integers $b\geq 2$.
\vskip 5pt

Theorem 2.2 will be proved in Section 6.
\vskip 5pt

For all integers $\alpha,\beta\geq 0$, define
$$
\displaylines{
\hfill F_{\alpha,\beta}(z)
= {1\over z^{2^\alpha}} \sum_{n=0}^\infty {z^{2^{n+\alpha}} \over 1+z^{2^{n+\beta}}}
=\sum_{n,j\geq 0}^{\infty} (-1)^jz^{(j 2^{\beta-\alpha}+1)2^{n+\alpha}-2^\alpha},\hfill \cr
\hfill G_{\alpha,\beta}(z)
= {1\over z^{2^\alpha}} \sum_{n=0}^\infty {z^{2^{n+\alpha}} \over 1-z^{2^{n+\beta}}}
=\sum_{n,j\geq 0}^{\infty} z^{(j 2^{\beta-\alpha}+1)2^{n+\alpha}-2^\alpha}.\hfill \cr
}
$$
The radius of convergence of $F_{\alpha,\beta}$ (resp. $G_{\alpha,\beta}$) is at least equal to
$1$. Moreover if $\beta=\alpha+1$, then $G_{\alpha,\beta}(z)$ is a rational function, since we have
$$
\leqalignno{
G_{\alpha,\alpha+1}(z)
&= {1\over z^{2^\alpha}} \sum_{n=0}^{\infty} {z^{2^{n+\alpha}}\over 1-z^{2^{n+\alpha+1}}}\cr
&= {1\over z^{2^\alpha}} \sum_{n=0}^{\infty}z^{2^{n+\alpha}}\sum_{j=0}^{\infty}z^{j 2^{n+\alpha+1}}\cr
            &= {1\over z^{2^\alpha}} \sum_{n=0}^{\infty}\sum_{j=0}^{\infty}(z^{2^\alpha})^{(2j+1)2^{n}}\cr
            &= {1\over z^{2^\alpha}} \sum_{m=0}^{\infty}(z^{2^\alpha})^{m}\cr
            &={1\over z^{2^\alpha}(1-z^{2^\alpha})}.\cr
}
$$
For $\beta\neq\alpha+1$, we have the following result.

\proclaim Theorem 2.3.
Let $\alpha,\beta\geq 0$ be integers such that $\beta\neq \alpha+1$.
Let $b\geq 2$ be an integer. Then both $F_{\alpha,\beta}(1/b)$ and $G_{\alpha,\beta}(1/b)$ 
are transcendental, and their irrationality exponent
are equal to $2$.

\noindent {\bf Remark.} The case $\alpha=\beta=0$  
implies that both ${\cal F}(1/b)$ and ${\cal G}(1/b)$ are transcendental for all integers $b\geq 2$,
and also the result obtained by Coons [Co13], namely that 
$\mu({\cal F}(1/b))=\mu({\cal G}(1/b))=2$, for all integers $b\geq 2$.
The case $\alpha=0$ and $\beta=2$ shows that for all integers $b\geq 2$,   
the regular paperfolding numbers $f_{RPF}(1/b)$
are transcendental and their irrationality exponents are equal to $2$.
The latter was conjectured by Coons and Vrbik~[CV12]
and has recently been established by Guo, Wu and Wen [GWW14].
\vskip 5pt

Recall that Stern's sequence $(a_n)_{n\geq 0}$
and its twisted version $(b_n)_{n\geq 0}$ are defined respectively by
(see [BV13, Ba10, St58])
$$
\cases{
a_0=0,\ a_1=1, \cr
a_{2n}=a_n,\ a_{2n+1}=a_n+a_{n+1}, \ (n\geq 1), \cr
}
$$
and
$$
\cases{
b_0=0,\ b_1=1,\cr
b_{2n}=-b_n, \   b_{2n+1}=-(b_n+b_{n+1}), \ (n\geq 1).
}
$$
Put $S(z)=\sum\limits_{n=0}^{\infty} a_{n+1}z^n$
and $T(z)=\sum\limits_{n=0}^{\infty} b_{n+1}z^n$.
Then $S$ and $T$  converge inside the unit disk, since $|a_n|\leq n$
and $|b_n|\leq n$ for all integers $n\geq 0$.
Recently, Bundschuh and V\"a\"an\"anen [BV13] proved
that $\mu(S(1/b))\leq 2.929$ and $\mu(T(1/b))\leq 3.555$
for all integers $b\geq 2$. Our next result gives the exact irrationality exponent
of the Stern number and also that of the twisted Stern number, and it will be proved in Section 7.

\proclaim Theorem 2.4.
For all integers $b\geq 2$, both $S(1/b)$ and $T(1/b)$ are transcendental and
their irrationality
exponents are equal to $2$.

The following theorem will be proved in Section 6.

\proclaim Theorem 2.5.  
Let $f(z)\in \setZ[[z]]$ be a power series defined by
$$
f(z)=\prod_{n=0}^{\infty}{C(z^{3^n})\over D(z^{3^n})},\leqno{(2.3)}
$$
with $D(z), C(z)\in \setZ[z]$ such that
$C(0)=D(0)= 1$.
Let $b\geq 2$ be an integer such that $C({1\over {b^{3^m}} })D({1\over {b^{3^m}} })\neq 0$
for all integers $m\geq 0$. If $f(z)\, (\rm{mod}\, 3)$ is not a rational function,
then $f(1/b)$ is transcendental and its irrationality exponent
is equal to $2$.

\noindent {\bf Remark.}
Taking $C(z)=1+z^2$ and $D(z)=1$ in Theorem 2.5 and
using the fact that the Cantor sequence on $\{0,1\}$ is not ultimately periodic,
we obtain that $f_C(1/b)$ is transcendental for all integers $b\geq 2$,
where the function $f_C$ is defined in Section 1. 
We also recover the result of [WW14] about Cantor real numbers, namely that
$\mu(f_C(1/b))=2$ for all integers $b\geq 2$. 

\vskip 5pt

For additional results, see Theorems 4.2, 6.1, 7.1, 7.2, and  Corollary 6.2.

\section{3. Hankel determinants and Pad\'e approximation} 

In this section we summarize several basic facts on
Pad\'e approximation. For more details, we refer the reader for example
to [Br80, BG96].

\medskip
Let $\setF$ be a field and $z$ be an indeterminate over $\setF$.
For any sequence ${\bf c}=(c_m)_{m\geq 0}$ of elements in $\setF$,
we put $f=f(z)=\sum\limits_{m=0}^{+\infty}c_m z^m$,
and call it the generating function of ${\bf c}$.
For all integers $n\geq 1$ and $k\geq 0$,
the Hankel determinant
of the power series $f$ (or of the sequence $\bf c$)  is defined by
$$
H_n^{(k)} (f) := \left|
\matrix{ c_k & c_{k+1} & \ldots & c_{k+n-1} \cr
c_{k+1} & c_{k+2} & \ldots & c_{k+n} \cr
\ \vdots \hfill & \ \vdots \hfill & \ddots &
\ \vdots \hfill \cr
c_{k+n-1} & c_{k+n} & \ldots & c_{k+2n-2} \cr} \right| \in \setF.
\leqno{(3.1)}
$$
By convention, we put $H_0^{(k)} (f)=1$, for all integers $k \ge 0$.
For all integers $n \ge 0$, write $H_n(f):=H_n^{(0)}(f)$. The sequence $H(f):=(H_n(f))_{n\geq 0}$
is called the {\it sequence of the Hankel determinants} of $f$.

Let $p$ and $q$ be nonnegative integers. By definition,
the Pad\'e approximant $[p/q]_f (z)$ to $f$ is the rational
fraction $P(z)/Q(z)$ in $\setF[[z]]$ such that
$$
\deg(P) \le p, \ \deg(Q) \le q,
\ \hbox{and $f(z)-{{P(z)}\over{Q(z)}}={\cal O}(z^{p+q+1})$}.
$$
The pair $(P,Q)$ has no reason to be unique, but the
fraction $P(z)/Q(z)$ is unique. Moreover if
we assume that $P$ and $Q$ are coprime, then $Q(0)\neq 0$.

If there exists an integer $k\geq 1$ such that $H_k (f)$ is nonzero,
then we know that the Pad\'e approximant
$[k-1 / k]_f (z)$ exists and we have
$$
f(z) - [k-1 / k]_f (z) = {H_{k+1} (f) \over H_k (f)} \, z^{2k}
+ {\cal O}(z^{2k+1}).  \leqno (3.2)
$$
This formula is of little help if $H_{k+1} (f) = 0$. But even in this case,
we still have the following fundamental result.

\proclaim Theorem 3.1.
With  the notation as above, suppose that there exist two integers $\ell,k$
such that $\ell > k\geq 1$ and
$H_{\ell}(f)H_k(f)\neq 0$. Then the Pad\'e approximant $[k-1 / k]_f (z)$ exists, and
there exist a nonzero element $h_k$ in~$\setF$
and an integer $k'$ such that $k\leq k'< \ell$ and
$$
f(z) - [k-1 / k]_f (z) = h_k \, z^{k+k'}+ {\cal O}(z^{k+k'+1}).  \leqno (3.3)
$$

\noindent {\bf Remark.} It seems to us that Theorem 3.1 is new.
An important point in its statement is
the non-vanishing of $h_k$.
\vskip 5pt

\noindent {\it Proof}.
Since $H_{\ell} (f)$ is nonzero, all the column
vectors in $H_{\ell} (f)$ are linearly independent,
in particular, the rank of the $\ell \times (k+1)$ matrix
$$
\pmatrix{ c_0 & c_1 & \ldots & c_{k} \cr
c_1    & c_2    & \ldots & c_{k+1} \cr
\vdots & \vdots & \ddots &
\vdots \cr
c_{\ell - 2} & c_{\ell - 1} & \ldots & c_{\ell + k - 2} \cr
c_{\ell-1} & c_{\ell} & \ldots & c_{\ell+ k - 1} \cr}
$$
is equal to $k+1$. By hypothesis, we also have $H_{k} (f)\neq 0$,
thus there exists a smallest integer $k'$
such that $k\leq k'<\ell$ and
$$
H_{k, k'} (f) := \left|
\matrix{
c_0     &c_1      &\ldots &c_{k}    \cr
c_1     &c_2      &\ldots &c_{k+1}  \cr
\vdots   &\vdots   &\ddots &\vdots   \cr
c_{k-1} &c_k      &\ldots &c_{2k-1} \cr
c_{k'}  &c_{k'+1} &\ldots &c_{k+k'} \cr} \right|\neq 0.
$$
Hence for all integers $j = k, \ldots, k'-1,$ we have $H_{k, j} (f)=0$. Define
$$
\leqalignno{Q^{[k-1/k]}(z)&:= \left|
\matrix{
c_0     &c_1     &\ldots &c_{k-1}  &c_{k} \cr
c_1     &c_2     &\ldots &c_k      &c_{k+1} \cr
\vdots  &\vdots  &\ddots &\vdots   &\vdots\cr
c_{k-1} &c_k     &\ldots &c_{2k-2} &c_{2k-1} \cr
z^k     &z^{k-1} &\ldots &z        &1 \cr}
\right|,\cr
P^{[k-1/k]}(z)&:= \left|
\matrix{
c_0     &c_1     &\ldots &c_{k-1}  &c_{k} \cr
c_1     &c_2     &\ldots &c_k      &c_{k+1} \cr
\vdots  &\vdots  &\ddots &\vdots   &\vdots\cr
c_{k-1} &c_k     &\ldots &c_{2k-2} &c_{2k-1} \cr
0   &c_0z^{k-1} &\ldots &\sum\limits_{i=0}^{k-2}c_iz^{i+1}       &\sum\limits_{i=0}^{k-1}c_iz^i\cr}
\right|.\cr
}
$$
So $\deg(P^{[k-1/k]})\leq k-1$, $\deg(Q^{[k-1/k]})\leq k$, and then (see [BG96, p.~6])
$$
\leqalignno{
&\quad Q^{[k-1/k]}(z)f(z)-P^{[k-1/k]}(z)\cr
&=Q^{[k-1/k]}(z)\Big(\sum\limits_{i=0}^{+\infty}c_iz^i\Big)-P^{[k-1/k]}(z)\cr
&= \left|
\matrix{
c_0     &c_1     &\ldots &c_{k-1}  &c_{k} \cr
c_1     &c_2     &\ldots &c_k      &c_{k+1} \cr
\vdots  &\vdots  &\ddots &\vdots   &\vdots\cr
c_{k-1} &c_k     &\ldots &c_{2k-2} &c_{2k-1} \cr
\sum\limits_{i=0}^{+\infty}c_iz^{i+k}     &\sum\limits_{i=0}^{+\infty}c_iz^{i+k-1} &\ldots
&\sum\limits_{i=0}^{+\infty}c_iz^{i+1}  &\sum\limits_{i=0}^{+\infty}c_iz^i \cr}
\right|-P^{[k-1/k]}(z)\cr
&=\left|
\matrix{
c_0     &c_1     &\ldots &c_{k-1}  &c_{k} \cr
c_1     &c_2     &\ldots &c_k      &c_{k+1} \cr
\vdots  &\vdots  &\ddots &\vdots   &\vdots\cr
c_{k-1} &c_k     &\ldots &c_{2k-2} &c_{2k-1} \cr
\sum\limits_{i=k}^{+\infty}c_iz^{i+k}     &\sum\limits_{i=k+1}^{+\infty}c_iz^{i+k-1} &\ldots
&\sum\limits_{i=2k-1}^{+\infty}c_iz^{i+1}  &\sum\limits_{i=2k}^{+\infty}c_iz^i \cr}
\right|\cr
&=\sum\limits_{i=1}^{+\infty}z^{2k+i-1}\left|
\matrix{
c_0     &c_1     &\ldots &c_{k-1}  &c_{k} \cr
c_1     &c_2     &\ldots &c_k      &c_{k+1} \cr
\vdots  &\vdots  &\ddots &\vdots   &\vdots\cr
c_{k-1} &c_k     &\ldots &c_{2k-2} &c_{2k-1} \cr
c_{k+i-1}     &c_{k+i} &\ldots &c_{2k+i-2}  &c_{2k+i-1} \cr}
\right|, \cr
}
$$
where in the first determinant, we have subtracted $z^k$ times the first row from the last one,
$z^{k+1}$ times the second row from the last one, etc.,
up to $z^{2k-1}$ times the penultimate row from the last one,
and then we arrive at the second determinant.

By the definition of the integer $k'$, we obtain
$$
Q^{[k-1/k]}(z)f(z)-P^{[k-1/k]}(z)=H_{k,k'}(f)z^{k+k'}+{\cal O}(z^{k+k'+1}).
$$
Note that $Q^{[k-1/k]}(0)=H_k(f)\neq 0$, thus we have
\goodbreak
$$
\leqalignno{f(z)-{{P^{[k-1/k]}(z)}\over {Q^{[k-1/k]}(z)}}
&={{H_{k,k'}(f)}\over {Q^{[k-1/k]}(z)}}z^{k+k'}+{\cal O}(z^{k+k'+1})\cr
&={{H_{k,k'}(f)}\over {H_k(f)}}z^{k+k'}+{\cal O}(z^{k+k'+1}).}
$$
Finally it suffices to put
$$
[k-1/k]_f(z):={{P^{[k-1/k]}(z)}\over {Q^{[k-1/k]}(z)}},\ h_k:={{H_{k,k'}(f)}\over {H_k(f)}}\neq 0,
$$
and we obtain at once the desired result. \qed
\medskip

To conclude this section, we recall some properties of rational functions in $\setZ[[z]]$, which are
related to Hankel determinants.
Let $(c_m)_{m\geq 0}$ be an integer sequence such that the power series
$f(z)=\sum\limits_{m=0}^{+\infty}c_mz^m$ converges inside the unit disk.
By Fatou's theorem (see [Fa06]), we know that the power series $f(z)$ is either rational
or transcendental over $\setQ(z)$. Moreover, by Kronecker's theorem
(see for example [Sa63, p.~5]),
we know also that the power series $f(z)$ is rational if and only
if there exists an integer $n_0\geq 0$
such that $H_n(f)=0$ for all integers $n$ larger than $n_0$. Equivalently, $f(z)$ is not rational
if and only if there exists an increasing sequence of positive integers
$(n_i)_{i\geq 0}$  such that $H_{n_i}(f)\not=0$, for all integers $i\geq 0$.

Finally we point out that since the power series $f(z)$ has only integer coefficients, thus it is 
transcendental over $\setQ(z)$
if and only if it is transcendental over $\setC(z)$ (see for example [SW88]).

\section{4. Irrationality exponent with Hankel determinants} 

In this section, we compute with Hankel determinants
the irrationality exponent of transcendental numbers, which
are special values at the inverse of integers $\geq 2$ of power
series satisfying a special type of functional equation.
\vskip 10pt

\proclaim Theorem 4.1.
Let $d\geq 2$ be an integer, and $(c_j)_{j \geq 0}$
be an integer sequence such that
$f(z)=\sum\limits_{j=0}^{+\infty}c_j z^j$ converges inside the unit disk.
Suppose that there exist integer polynomials
$A(z)$, $B(z)$, $C(z)$ and $D(z)$ such that
$$
f(z)={A(z)\over B(z)} + {C(z)\over D(z)} f(z^d).\leqno{(4.1)}
$$
Let $b\geq 2$ be an integer such that $B({1\over {b^{d^m}}})C({1\over {b^{d^m}}})D({1\over {b^{d^m}}})\neq 0$,
for all integers $m\geq 0$.
If there exists an increasing sequence of positive integers
$(n_i)_{i\geq 0}$  such that $H_{n_i}(f)\not=0$ for all integers $i\geq 0$ and
$\limsup\limits_{i\rightarrow +\infty} {n_{i+1}\over n_i} =\rho $,
then $f(1/b)$ is transcendental, and we have
$$
\mu\Big(f\Big({1\over b}\Big)\Big)\leq (1+\rho)\min\{\rho^2, d\}.
$$
In particular, the irrationality exponent of $f(1/b)$ is equal to $2$ if $\rho=1$.

\noindent {\it Proof}.
From the equation (4.1), we deduce immediately that 
$$
\left(\matrix{1\cr
f(z)}\right)=\left(\matrix{1&0\cr
{A(z)\over B(z)}&{C(z)\over D(z)}}\right)\left(\matrix{1\cr
f(z^d)}\right).
$$
Since $B({1\over {b^{d^m}}})C({1\over {b^{d^m}}})D({1\over {b^{d^m}}})\neq 0$ for all integers $m\geq 0$,
then by a result due to Nishioka (see [Ni90, Corollary 2]), we obtain
$$
{\rm tr}.\deg_{\setQsmall}\setQ\big(1, f(1/b)\big)
={\rm tr}.\deg_{\setCsmall(z)}\setC(z)\big(1, f(z)\big).
$$
Now that there exists an increasing sequence of positive integers
$(n_i)_{i\geq 0}$ such that $H_{n_i}(f)\not=0$ for all integers $i\geq 0$,
the power series $f(z)$ is not a rational function.  
Thus, it is transcendental over $\setC(z)$ by Fatou's theorem,  
hence $f(1/b)$ is transcendental.

By iteration of Formula (4.1), we have, for all integers $m\geq 1$,   
$$
f(z)={{A_m}(z)\over {B_m(z)}}+{{C_m}(z)\over {D_m(z)}}f(z^{d^m}),\leqno{(4.2)}
$$
where $C_m(z)=\prod\limits_{j=0}^{m-1}C(z^{d^j})$,
$D_m(z)=\prod\limits_{j=0}^{m-1}D(z^{d^j})$, and
$$
B_m(z)=D_{m-1}(z)\prod\limits_{j=0}^{m-1}B(z^{d^j}),\
A_m(z)=\sum\limits_{j=0}^{m-1}C_{j}(z)A(z^{d^j})\cdot {{B_m(z)}\over {D_{j}(z)B(z^{d^j})}},
$$
where we have put $C_0(z) = D_0(z) = 1$. 

Put $\alpha=\deg(A(z))$, $\beta=\deg(B(z))$, $\gamma=\deg(C(z))$, $\delta=\deg(D(z))$. Then,
$$
\leqalignno{
\deg(C_m(z))&=\sum\limits_{j=0}^{m-1}\deg(C(z^{d^j}))
=\sum\limits_{j=0}^{m-1}\gamma d^j ={{\gamma(d^m-1)}\over {d-1}}\leq \gamma d^m,\cr
\deg(D_m(z))&=\sum\limits_{j=0}^{m-1}\deg(D(z^{d^j}))
=\sum\limits_{j=0}^{m-1}\delta d^j ={{\delta(d^m-1)}\over {d-1}}\leq \delta d^m,\cr
\deg(B_m(z))&=\deg(D_{m-1}(z))+\sum\limits_{j=0}^{m-1}\deg(B(z^{d^j}))\cr
            &={{\delta(d^{m-1}-1)}\over {d-1}}+{{\beta(d^m-1)}\over {d-1}}\leq (\beta+\delta)d^m,\cr
\deg(A_m(z))
&\leq \max_{0\leq j\leq m-1}\Big(\deg(C_j(z))+\deg(A(z^{d^j}))+\deg(B_m(z)) \Big)\cr 
            &\leq \max_{0\leq j\leq m-1}\Big({{\gamma(d^j-1)}\over {d-1}}+\alpha d^j+{{\delta(d^{m-1}-1)}\over {d-1}} +{{\beta(d^m-1)}\over {d-1}} \Big)\cr
            &\leq (\alpha+\beta+\gamma+\delta)d^m.    
}
$$
Let $i\geq 0$ be an integer.
As in the proof of Theorem 3.1,
we denote by $n_i'$ the smallest integer such
that $n_i\leq n_i'<n_{i+1}$ and $H_{n_i,n_i'}(f)\neq 0$.
Then we can find $h_i\in\setQ\setminus \{0\}$, and
$P_i(z),Q_i(z)\in \setZ[z]$ with $\deg(P_i(z))\leq n_i-1$, $\deg(Q_i(z))\leq n_i$,
and $Q_i(0)\neq 0$ such that
$$
f(z)-{{P_i(z)}\over{Q_i(z)}}=h_iz^{n_i+n_i'}+{\cal O}(z^{n_i+n_i'+1})
=h_iz^{n_i+n_i'}\big(1+{\cal O}(z)\big).
$$
Thus, for all integers $m\geq 1$, we obtain
$$
f(z^{d^m})-{{P_i(z^{d^m})}\over{Q_i(z^{d^m})}}=h_iz^{(n_i+n_i')d^m}\big(1+{\cal O}(z^{d^m})\big).
$$
Combined with Formula (4.2), this gives
$$
 f(z)-{{A_m}(z)\over {B_m(z)}}-{{C_m}(z)\over {D_m(z)}}\cdot {{P_i(z^{d^m})}\over{Q_i(z^{d^m})}}
=h_iz^{(n_i+n_i')d^m}{{C_m}(z)\over {D_m(z)}}\big(1+{\cal O}(z^{d^m})\big).
$$
To simplify the notation, we define
$$
\leqalignno{
P_{i,m}(z)&=A_m(z)D_m(z)Q_i(z^{d^m})-B_m(z)C_m(z)P_i(z^{d^m}),\cr
Q_{i,m}(z)&=B_m(z)D_m(z)Q_i(z^{d^m}).
}
$$
Since $B(z)C(z)D(z)\neq 0$, then we can write
$$
B(z)=b_\kappa z^{\kappa} (1+z \tilde{B}(z)),\
C(z)=c_\eta z^{\eta} (1+z \tilde{C}(z)),\
D(z)=d_{\iota}z^{\iota}(1+z\tilde{D}(z))
$$
with $\kappa,\eta, \iota\geq 0$ integers, $b_\kappa, c_{\eta},d_{\iota}\in\setZ$,
and $\tilde{B}(z), \tilde{C}(z), \tilde{D}(z)\in\setQ[z]$.
Note that $\tilde{C}(z), \tilde{D}(z)$ are bounded on the unit disk, thus both
$\sum\limits_{j=0}^{\infty}{1\over b^{d^j}}|\tilde{C}({1\over b^{d^j}})|$ and
$\sum\limits_{j=0}^{\infty}{1\over b^{d^j}}|\tilde{D}({1\over b^{d^j}})|$ converge.
Note also that $C({1\over b^{d^m}})D({1\over b^{d^m}})\neq 0$ for all integers $m\geq 0$,
thus the following two limits
$$
\leqalignno{
\sigma &=\lim_{m\rightarrow +\infty}
{{C_m({1\over b})}\over {c_\eta^m b^{-{{\eta(d^m-1)}\over {d-1}}}}}
=\prod_{j=0}^{\infty}\Big(1+{1\over b^{d^j}}\tilde{C}\big({1\over b^{d^j}}\big)\Big),\cr
\tau &=\lim_{m\rightarrow +\infty}{{D_m({1\over b})}\over {d_\iota^m b^{-{{\iota(d^m-1)}\over {d-1}}}}}
=\prod_{j=0}^{\infty}\Big(1+{1\over b^{d^j}}\tilde{D}\big({1\over b^{d^j}}\big)\Big)\cr
}
$$
do exist and are different from zero.
Hence, for $m$ tending to $+\infty$, we have
$$
\eqalign{
f\Big({1\over b}\Big)-{{P_{i,m}({1\over b})}\over {Q_{i,m}({1\over b})}}
&={h_i\over b^{(n_i+n_i')d^m}}
{{C_m}({1\over b})\over {D_m({1\over b})}}\Big(1+{\cal O}\big({1\over b^{d^m}}\big)\Big)\cr
&\sim {{h_i\sigma}\over{\tau}}\Big({{c_{\eta}}\over {d_{\iota}}}\Big)^m
{1\over b^{(n_i+n_i')d^m+{{(\eta+\iota)(d^m-1)}\over {d-1}}}}.\cr}   \leqno{(4.3)}
$$
Moreover, we also have
$$
\leqalignno{\deg(P_{i,m}(z))&\leq \max\Big(\deg(A_m(z))+\deg(D_m(z))+\deg(Q_i(z^{d^m})),\cr
&\hskip 10pt \deg(B_m(z))+\deg(C_m(z))+\deg(P_i(z^{d^m}))\Big)\cr
&\leq \max\Big((\alpha+\beta+\gamma+2\delta+n_i)d^m,
(\beta+\delta+\gamma+n_i-1)d^m))\Big)\cr
&\leq (\alpha+\beta+\gamma+2\delta+n_i)d^m,\cr
\deg(Q_{i,m}(z))&\leq \deg(B_m(z))+\deg(D_m(z))+\deg(Q_i(z^{d^m}))\cr
&\leq (\beta+2\delta+n_i)d^m.
}
$$
Put $e_i=\alpha+\beta+\gamma+2\delta+n_i$.
Recall that, by assumption, we have
$B({1\over {b^{d^m}}})D({1\over {b^{d^m}}})\neq 0$
for all integers $m\geq 0$.
Since $Q_i(0)\neq 0$,
we can find two constants $\alpha_{1,i},\alpha_{2,i}>0$
(which depend only on $i$)
such that for all integers $m\geq 0$, we have
$$
\alpha_{1,i}\leq
b^{\kappa d^m}\Big| B\Big({1\over {b^{d^m}}}\Big) \Big|,
b^{\iota d^m}\Big| D\Big({1\over {b^{d^m}}}\Big) \Big|,
\Big| Q_i\Big({1\over {b^{d^m}}}\Big) \Big|
\leq \alpha_{2,i}.
$$
Set $q_{i,m}=b^{e_id^m}|Q_{i,m}({1\over b})|$, and $p_{i,m}=b^{e_id^m}P_{i,m}({1\over b}){\rm sgn}(Q_{i,m}({1\over b}))$.
Then $q_{i,m},p_{i,m}$ are integers, and for all integers $m\geq 1$, we have
$$
\alpha_{1,i}^{3m}b^{e_id^m-g_m}
\leq q_{i,m}
\leq \alpha_{2,i}^{3m}b^{e_id^m-g_m}, \leqno{(4.4)}
$$
with $g_m={{(\kappa+\iota)(d^m-1)+\iota(d^{m-1}-1)}\over {d-1}}$,
from which we deduce immediately
$$
{{\alpha_{1,i}^{3(m+1)}}\over {\alpha_{2,i}^{3m}}}b^{e_id^m(d-1)-(d\kappa+d\iota+\iota)d^{m-1}}q_{i,m}
\leq q_{i,m+1}\leq {{\alpha_{2,i}^{3(m+1)}}\over {\alpha_{1,i}^{3md}}}q_{i,m}^d.\leqno{(4.5)}
$$
Let $\varepsilon$ be a sufficiently small positive real number.  
Since $\lim\limits_{i\rightarrow +\infty}n_i=+\infty$,
there exists an integer $N_1>1$ (independent of $m$) 
such that for $i>N_1$ and $m \ge 2$, we have    
$$\cases{
\displaystyle (1-\eps)n_id^m\leq e_id^m-g_m\leq (1+\eps)n_id^m, \cr
\displaystyle e_id^m(d-1)-(d\kappa+d\iota+\iota)d^{m-1}>(1-\eps)n_id^m.}
$$ 
Then it follows from Formulas (4.4) and (4.5)
that there exists an integer $N_{1,i}>1$
such that for all integers $m\geq N_{1,i}$, we have
$$
\displaylines{
\rlap{(4.6)}\hfill b^{n_id^m(1-2\varepsilon)}\leq q_{i,m}\leq b^{n_id^m(1+2\varepsilon)},\hfill\cr
\rlap{(4.7)}\hfill q_{i,m}<q_{i,m+1}\leq  q_{i,m}^{d(1+\varepsilon)}.\hfill
}
$$
Similarly it follows from Formula (4.3) that there exists
an integer $N_{2,i}>N_{1,i}$
such that for all integers $m\geq N_{2,i}$, we have
$$
{1 \over {b^{(n_i+n_i'+\eta+\iota)(1+\varepsilon)d^m}}}
\leq \Big|f\Big({1\over b}\Big)-{{p_{i,m}}\over {q_{i,m}}}\Big|
\leq {1 \over {b^{(n_i+n_i')(1-\varepsilon)d^m}}} \leqno{(4.8)}
$$
and by Formula (4.6), we obtain also
$$
{1 \over {q_{i,m}^{(n_i+n_i'+\eta+\iota)(1+ 4 \varepsilon)/n_i }}}
\leq \Big|f\Big({1\over b}\Big)-{{p_{i,m}}\over {q_{i,m}}}\Big|
\leq {1 \over {q_{i,m}^{(n_i+n_i')(1- 4 \varepsilon)/ n_i}}}. \leqno{(4.9)}    
$$
By hypothesis, we have $\limsup\limits_{i\rightarrow\infty}{{n_{i+1}}\over {n_i}}=\rho $,
then we can find an integer $i_0> N_1$
such that for all integers $i\geq i_0$, we have ${{n_{i+1}}\over {n_i}}<\rho+\varepsilon$
and
$$
\cases{\displaystyle
{(n_i+n_i'+\eta+\iota)(1+ 4 \varepsilon) \over n_i } \le (1+\rho) (1 + 6 \eps), \cr    
\noalign{\smallskip}
\displaystyle {(n_i+n_i')(1- 4 \varepsilon) \over n_i} \ge 2 (1 - 6 \eps), \cr} \leqno{(4.10)}  
$$
from which we deduce at once
$$
{1 \over {q_{i,m}^{(1+\rho)(1+ 6 \varepsilon) }}}
\leq \Big|f\Big({1\over b}\Big)-{{p_{i,m}}\over {q_{i,m}}}\Big|
\leq {1 \over {q_{i,m}^{2(1- 6 \varepsilon)}}}. \leqno{(4.11)}
$$
Applying Lemma 4.1 from [AR09, p.~668] with (4.7) and (4.11), we obtain
$$
\mu\Big(f\Big({1\over b}\Big)\Big)\leq {(1+\rho) (1 + 6 \eps)\over 2 (1 - 6 \eps)-1}d(1+\varepsilon).
$$
Since $\varepsilon$ is positive and can be chosen arbitrarily small, we get
$$
\mu\Big(f\Big({1\over b}\Big)\Big)
\leq (1+\rho)d.   \leqno (4.12)
$$

Fix $\ell >1$ an integer such that $d^{\ell-1}>n_{i_0}$. Let ${\cal A}_{\ell}$
be the set of integers $i>i_0$ such that $n_i\in [d^{\ell-1}, d^{\ell}-1]$.
Assume that ${\cal A}_{\ell}$ is non-empty (it could be empty when $\rho$ is large,   
but it is certainly non-empty for infinitely many~$\ell$),
and denote its elements as $n_{i_1}<n_{i_2}<\cdots<n_{i_t}$. Then $t\geq 1$,
$n_{i_j}=n_{i_1+j-1}\ (1\leq j\leq t)$, $n_{i_1}<(\rho+\eps)n_{i_1-1}<(\rho+\eps)(d^{\ell-1}-1)$, and
$d^{\ell}\leq n_{i_t+1}<(\rho+\varepsilon)n_{i_t}$. Put
$$
M_{\ell}=\max_{1\leq i\leq i_t}N_{2,i}.
$$
Arrange the integers $q_{i_l,m}$ ($1\leq l\leq t$ and $m\geq M_{\ell}$) as an increasing
sequence, which we denote by $(r_{\ell,j})_{j\geq 0}$.

Fix $j\geq 0$, and write
$r_{\ell,j} = q_{i_l,m}$ with $1\leq l\leq t$. By (4.6), we have
$$
b^{n_{i_l}  d^m(1-2\varepsilon)}\leq q_{i_l,m}\leq b^{n_{i_l} d^m(1+2\varepsilon)}.
$$

We distinguish below two cases:
\smallskip

{\bf Case I:} $n_{i_{t}} > n_{i_l} (1 + 2\varepsilon) / (1 - 2\varepsilon)$.
Then $i_{t}>i_l$, and thus there exists a smallest integer
$v$ such that $l<v\leq t$ such that
$$
n_{i_v} > n_{i_l} (1 + 2\varepsilon) / (1 - 2\varepsilon).
$$
Consequently we have $q_{i_l,m} < q_{i_v, m}$ and
$$
{\log q_{i_v,m}\over \log q_{i_l,m}}
\le  {n_{i_v}  (1 + 2\varepsilon)\over n_{i_l} (1 - 2\varepsilon)}.
$$
By the minimality $v$, we have
$$
n_{i_v} <   (\rho + \varepsilon) n_{i_{v-1}} \le
(\rho + \varepsilon)  n_{i_l} (1 + 2\varepsilon) / (1 - 2\varepsilon),  
$$
from which we deduce directly
$$
1<{\log r_{\ell,j+1}\over \log r_{\ell,j}}\leq {\log q_{i_v,m}\over \log q_{i_l,m}}
 <  {( \rho  + \varepsilon)   (1 + 2\varepsilon)^2\over (1 - 2\varepsilon)^2}.
$$

{\bf Case II:} $n_{i_t}  \le n_{i_l} (1 + 2\varepsilon) / (1 - 2\varepsilon)$.
Since $n_{i_t}<d^{\ell}\leq d  n_{i_1}$, we have
$$
n_{i_t} {1 + 2\eps \over 1 - 2 \eps} < d n_{i_1},
$$
for all $\eps>0$ small enough. Then we get
$$
{\log q_{i_1,m+1}\over \log q_{i_l,m}}
\geq {n_{i_1} d(1-2\eps) \over n_{i_l}(1+2\varepsilon)}>{n_{i_t}\over n_{i_l}}\geq 1.
$$
Moreover, from $n_{i_t}  \le n_{i_l} (1 + 2\varepsilon) / (1 - 2\varepsilon)$,
we obtain also
$$
{\log q_{i_1,m+1}\over \log q_{i_l,m}}
\leq { d n_{i_1} (1 + 2\eps) \over n_{i_l} (1 - 2 \eps)}
\le {d n_{i_1} (1 + 2\eps)^2 \over n_{i_t} (1 - 2 \eps)^2}.
$$
Note that $n_{i_1}<(\rho+\eps)(d^{\ell-1}-1)$
and $n_{i_t}>{d^{\ell}\over \rho+\eps}$, hence
$$
{dn_{i_1}\over n_{i_t}}<(\rho+\eps)^2{d(d^{\ell-1}-1)\over d^{\ell}}<(\rho+\eps)^2,
$$
and then we obtain
$$
1<{\log r_{\ell,j+1}\over \log r_{\ell,j}}\leq {\log q_{i_1,m+1}\over \log q_{i_l,m}}
 <  {(\rho+\eps)^2  (1 + 2\varepsilon)^2\over (1 - 2\varepsilon)^2}.
$$

In conclusion, since $\rho\geq 1$,
we have established in both cases that
$$
1<{\log r_{\ell,j+1}\over \log r_{\ell,j}}
 <  {(\rho+\eps)^2  (1 + 2\varepsilon)^2\over (1 - 2\varepsilon)^2},
 \leqno{(4.13)}
$$
for all integers $j\geq 0$.

Once again applying Lemma 4.1 from [AR09, p.~668] with (4.11) and (4.13), we get
$$
\mu\Big(f\Big({1\over b}\Big)\Big)\leq {(1+\rho)(1 + 6 \eps)\over 2 (1 - 6 \eps)-1} \cdot
{(\rho+\eps)^2  (1 + 2\varepsilon)^2\over (1 - 2\varepsilon)^2}.
$$
Since $\varepsilon$ is positive and can be chosen arbitrarily small, we obtain
$$
\mu\Big(f\Big({1\over b}\Big)\Big)\leq (1+\rho)\rho^2.
$$
Combined with (4.12), this gives
$$
\mu\Big(f\Big({1\over b}\Big)\Big)\leq (1+\rho) \min\{\rho^2, d\},  
$$
as asserted. 
In particular, if $\rho=1$, then $f(1/b)\le 2$. But $f(1/b)$ is transcendental,
thus its irrationality exponent is equal to $2$.
\qed
\vskip 10pt

\noindent {\bf Remarks.}
(1) Note that Nishioka's result (quoted at the beginning of the proof of Theorem 4.1) 
may fail if we remove the condition that
$C({1\over b^{d^m}})\neq 0$ for all integers $m\geq 0$.
Consider the power series
$$
f(z)=\prod\limits_{n\geq 0}(1-2z^{2^n}).
$$
Then $f(z)=(1-2z)f(z^2)$, and $f$ is analytic inside
the unit disk. It is also a transcendental function
for it has infinitely many zeros. However $f(1/2)=0$.
For more detail on this example, see [Be94, p.~283].

(2) In the statement of Theorem 4.2, if we replace $1/b$ by $a/b$ 
with $a$ an integer satisfying $0<|a|<b$,
then the same proof yields that $f(a/b)$ is transcendental.
If we suppose further $0<|a|<\sqrt{b}$, then with slight modifications, we can show
the upper bound (see [Du14] for the case of Thue--Morse)  
$$
\mu\Big(f\Big({a\over b}\Big)\Big)\leq {{\log b-\log |a|}\over {\log b-2\log |a|}}(1+\rho)\min\{\rho^2,d\}.
$$
\vskip 5pt

We are now in position to establish Theorem 2.1.
\medskip
\noindent {\it Proof of Theorem 2.1}.
Without loss of generality, we can suppose that the polynomials
$A(z)$ and $B(z)$ are coprime. From the functional equation, we obtain
that ${A(z)\over B(z)}=f(z)-C(z)f(z^d)$ is analytic inside the unit disk,
so $B({1\over {b^{d^m}}})\neq 0$ for all integers $m\geq 0$ and
$b\geq 2$. Then by Theorem 4.1, the desired result holds. \qed
\vskip 10pt

We display another application of Theorem 4.1.

\proclaim Theorem 4.2.
Let $d\geq 2$ be an integer, and $(c_j)_{j \geq 0}$ be an integer sequence
taking only finitely many values.
Put $f(z)=\sum\limits_{j=0}^{+\infty}c_j z^j$.
Suppose that there exist integer polynomials $A(z)$, $B(z)$,
$C(z)$ and $D(z)$ such that
$$
f(z)={A(z)\over B(z)} + {C(z)\over D(z)} f(z^d). 
$$
Let $b\geq 2$ be an integer such that $C\big({1\over b^{d^m}}\big)\neq 0$
for all integers $m\geq 0$.
If there exists an increasing sequence of positive integers
$(n_i)_{i\geq 0}$  such that $H_{n_i}(f)\not=0$ for all integers $i\geq 0$ and
$\lim\limits_{i\rightarrow \infty} {n_{i+1}\over n_i} =1$,
then $f(1/b)$ is transcendental and its
irrationality exponent of is equal to $2$.

\noindent {\it Proof}.
Since the sequence $(c_j)_{j\geq 0}$ is bounded,
the function $f(z)$ converges inside the unit disk, and for all integers $b\geq 2$, 
we can find an integer $\ell>2$ such that $|c_j|<b^{d^{\ell}-1}$, for all integers $j\geq 0$.
Note also that $f(z)$ is not rational, for  there exists
an increasing sequence of positive integers
$(n_i)_{i\geq 0}$  such that $H_{n_i}(f)\not=0$ for $i \ge 0$.
As in the proof of Theorem 4.1, for any integer $\ell > 2$,
we can find $A_{\ell}(z),B_{\ell}(z),C_{\ell}(z),D_{\ell}(z)$ in $\setZ[z]$ such that
$$
f(z)={A_{\ell}(z)\over B_{\ell}(z)} + {C_{\ell}(z)\over D_{\ell}(z)} f(z^{d^{\ell}}).\leqno{(4.14)}
$$
Without loss of generality, we can also suppose that
$$
\gcd(A_{\ell}(z),B_{\ell}(z))=1,\ {\rm and}\ \gcd(C_{\ell}(z),D_{\ell}(z))=1.
$$
We argue by contradiction.
Suppose that there is an integer $m\geq 0$ such that
$B_{\ell}\big({1\over b^{d^m}}\big)D_{\ell}\big({1\over b^{d^m}}\big)=0$.
Then we can write
$$
B_{\ell}(z)=\Big(z-{1\over b^{d^m}}\Big)^sE(z),\ D_{\ell}(z)=\Big(z-{1\over b^{d^m}}\Big)^tF(z),
$$
where $E(z), F(z)\in \setQ[z]$ are not equal to zero at $z={1\over b^{d^m}}$,
and $s,t\geq 0$ are integers such that $\max \{s,t\} \geq 1$.

If $s>t$, then from Formula (4.14), we obtain
$$
\Big(z-{1\over b^{d^m}}\Big)^tf(z)- {C_{\ell}(z)\over F(z)} f(z^{d^{\ell}})
={A_{\ell}(z)\over (z-{1\over b^{d^m}})^{s-t}E(z)}.
$$
The left hand side
is regular at $z={1\over b^{d^m}}$, while the right side is not, giving us
the required contradiction.

If $s\leq t$, then from Formula (4.14), we have
$$
\Big(z-{1\over b^{d^m}}\Big)^tf(z)- {(z-{1\over b^{d^m}})^{t-s}A_{\ell}(z)\over E(z)}
={C_{\ell}(z)\over F(z)} f(z^{d^{\ell}}).
$$
Hence $f(1/{b^{d^{m+\ell}}})$ is a rational number.
But $(c_j)_{j\geq 0}$ is the sequence of coefficients of this rational number in 
its  base-$b^{d^{m+\ell}}$ expansion and it 
is bounded by $b^{d^{\ell}-1}$. Thus, the sequence $(c_j)_{j\geq 0}$
is ultimately periodic. This gives again a contradiction since $f(z)$ is not rational.

To conclude, it suffices to apply Theorem 4.1 to the equation (4.14).~\qed
\vskip 10pt

The above theorems have many applications, but they also have an inconvenient:
in general it is not at all easy to check the conditions about Hankel determinants,
and indeed it is often extremely technical to compute explicitly
Hankel determinants (see for example [APWW98, GWW14]). Later we shall compute
the irrationality exponent only with information on the functional equation
satisfied by the related power series. For this, we need recall some basic
results about $J$-fractions in the following section.


\section{5. Hankel continued fraction} 

For proving Theorem 2.2, we need the {\it grafting} technique,
which has been introduced in [H15a] for the Jacobi continued fraction,
and extended for the  Hankel continued fraction  in [H15b].

For all integers $\delta\geq 1$, a {\it super continued fraction} associated with $\delta$,
called {\it super $\delta$-fraction} for short, is defined to be a continued fraction
of the following form (see [H15b]):
$$
f(z)
={v_0 z^{k_0}\over {
		1+u_1(z)z-\displaystyle{\strut v_1 z^{k_0+k_1+\delta}\over{
				1+	u_2(z)z-\displaystyle{\strut v_2 z^{k_1+k_2+\delta}\over{
						1+	u_3(z)z-\displaystyle{ \ddots
	}}}}}}}\leqno{(5.1)}
$$
where $v_j\not=0$ are constants, $k_j$ are nonnegative integers and $u_j(z)$ are polynomials of
degree less than or equal to $k_{j-1}+\delta-2$. By convention, we set $\deg\, 0=-1$.
\medskip

A super $2$-fraction is called an {\it Hankel continued fraction}.
The following two results about  Hankel continued fractions are established in [H15b].

\proclaim Theorem 5.1.
(i) Each  Hankel continued fraction defines
a power series, and conversely, for each power series $f(z)$,
the  Hankel continued fraction  expansion of $f(z)$ exists
and is unique.
\smallskip
(ii) Let $f(z)$ be a power series such that its
Hankel continued fraction is given by (5.1) with $\delta=2$.
Then, for all integers $j\geq 0$, all non-vanishing Hankel determinants of $f(z)$ are given by
$$
H_{s_j}(f)
= (-1)^\epsilon v_0^{s_j} v_1^{s_j-s_1} v_2^{s_j-s_2} \cdots v_{j-1}^{s_j-s_{j-1}}, \leqno{(5.2)}
$$
where $\epsilon= \sum\limits_{i=0}^{j-1} {k_i(k_i+1)/2}$
and
$s_j=k_0+k_1+\cdots + k_{j-1}+j$.
\medskip

For any prime number $p$, let $\setF_p:=\setZ/p\setZ$ denote the finite field
with $p$ elements.

\proclaim Theorem 5.2.
Let $p$ be a prime number and $F(z)\in \setF_p[[z]]$ be a power series
satisfying the following quadratic equation
$$
A(z)+B(z)F(z)+ C(z)F(z)^2 = 0, \leqno{(5.3)}
$$
where $A(z), B(z), C(z)\in \setF_p[z] $ are three polynomials satisfying
one of the following four conditions:
\smallskip 
(i) $B(0)=1, \ C(0)=0,\ C(z)\not=0$;
\smallskip 
(ii) $B(0)=1,\ C(z)=0$;
\smallskip 
(iii) $A(0)=0,\ B(0)=1,\ C(0)\not=0$;
\smallskip 
(iv) $p \ge 3$, $B(z)=0,\  C(0)=1$, and there exist an integer $k \ge 0$, $a_k$~in~$\setF_p\setminus\{0\}$, and $\tilde{A} (z)$
in $\setF_p[z]$ such that
$A(z)= -(a_k z^k)^2 (1 + z \tilde{A} (z))$.
\smallskip\noindent
Then, the Hankel continued fraction expansion of $F(z)$ exists and
is ultimately periodic. Also, the sequence of the Hankel determinants of $F$
is ultimately periodic.

\section{6. Irrationality exponent without Hankel determinants} 
\medskip

In this section, based on the information of the functional equation satisfied
by the power series and applying the results of the
previous section, we shall present several results about irrationality exponents
without explicit conditions on Hankel determinants.

\proclaim Theorem 6.1.  
Let $f(z)\in \setZ[[z]]$ be a power series analytic in the unit disk
and such that
$$
A(z)+B(z)f(z)+ C(z)f(z^2) = 0, \leqno{(6.1)}
$$
where $A(z)$, $B(z)$ and $C(z)$ are integer polynomials
satisfying one of the following conditions:
\smallskip
(i) $B(0)\equiv 1,\ C(0)\equiv 0\quad ({\rm mod}\, 2)$,
\smallskip
(ii) $A(0)\equiv 0,\  B(0)\equiv 1,\ C(0)\not\equiv 0\quad ({\rm mod}\, 2)$.
\smallskip\noindent
Let $b\geq 2$ be an integer such that $B({1\over {b^{2^m}} })C({1\over {b^{2^m}} })\neq 0$ for all integers $m\geq 0$.
If $f(z)\,(\rm{mod}\, 2)$ is not a rational function,
then $f(1/b)$ is transcendental and its irrationality exponent is equal to $2$.

\noindent {\it Proof}.
Put $F(z)=f(z)\, ({\rm mod}\, 2)\in \setF_2[[z]]$. By Formula (6.1), we obtain
$$
A(z)+B(z)F(z)+ C(z)F(z)^2 = 0.
$$
By Theorem 5.2 (with conditions (i) and (iii), respectively)
the sequence $H(F)$ of Hankel determinants is ultimatly periodic
over the field~$\setF_2$.
Since $F(z)$ is not a rational function in $\setF_2[[z]]$, there exists
an increasing sequence of positive integers
$(n_i)_{i\geq 0}$  such that
$H_{n_i}(F)\not=0$ for all integers $i\geq 0$ and
$\lim\limits_{i\rightarrow \infty} {n_{i+1}\over n_i} =1$.
Let $b\geq 2$ be an integer such that $B({1\over {b^{2^m}}})C({1\over {b^{2^m}}})\neq 0$ for all integers $m\geq 0$.
Then it follows from Theorem 4.1 that $f(1/b)$ is transcendental and its irrationality exponent is equal to $2$. \qed
\bigskip

\noindent {\it Proof of Theorem 2.5}.
Directly from the definition and the fact that $C(0)=D(0)=1$, we obtain
that $f(z)$ converges in the unit disk, its coefficients in power series expansion
are integers, and $f(z)={{C(z)}\over {D(z)}}f(z^3)$.
Over the field $\setF_3$, the power series $F(z)=f(z)\, (\rm{mod}\, 3)$ satisfies the
quadratic equation
$-D(z) + C(z)F(z)^2 =0$.
So by Theorem 5.2 (iv), the sequence $H(F)$ of Hankel determinants is ultimately periodic
over the field $\setF_3$.
Since $F(z)$ is not a rational function in $\setF_3[[z]]$, there exists
an increasing sequence of positive integers
$(n_i)_{i\geq 0}$  such that
$H_{n_i}(F)\not=0$ for all integers $i\geq 0$ and
$\lim\limits_{i\rightarrow \infty} {n_{i+1}\over n_i} =1$.
Let $b\geq 2$ be an integer such that $C({1\over {b^{3^m}}})D({1\over {b^{3^m}}})\neq 0$ for all integers $m\geq 0$.
It follows from Theorem 4.1 that $f(1/b)$ is transcendental and its irrationality exponent
is equal to $2$. \qed

\medskip

Letting $C(z)=1-z$ (resp. $C(z)=1\pm z-z^2$) and $D(z)=1$ in Theorem~2.5,
we obtain at once the following corollary. The underlying Hankel determinants are
evaluated in [H15a].

\proclaim Corollary 6.2.
For all integers $b\geq 2$, both
$$
\prod\limits_{k\geq 0} (1-b^{-3^k})
\quad \hbox{and} \quad
\prod\limits_{k\geq 0}(1\pm b^{-3^k} - b^{-2\cdot3^k})
$$
are transcendental and their irrationality exponents are equal to $2$.

\medskip

We are now in position to establish Theorem 2.2.
\vskip 5pt

\noindent {\it Proof of Theorem 2.2}.
From Formula  (2.2) and $D(0)=1$, we obtain directly that the power series $f(z)$ converges
in the unit disk,  its coefficients in power series expansion are integers, and
 $f(z)=\bigl(1+uz+2z^2{{C(z)}\over {D(z)}}\bigr)f(z^2)$,
$$
{{1}\over {f(z)}}= { {1-uz+\big(-2 C(0)+u^2-u\big)z^2}} + \cdots
$$
Since $u(u-1)$ is even, we can define $g(z)\in\setZ[[z]]$ by
$$
f(z)={1\over 1-uz+2z^2 g(z)}.\leqno{(6.2)}
$$
By Theorem 5.1 (ii) (or Lemma 2.2 in [H15b]),
the Hankel determinants of $f$ and those of $g$ are tightly related by
$$
H_n(f)=(-2)^{n-1}H_{n-1}(g).\leqno{(6.3)}
$$
By the functional equation satisfied by $f(z)$ and Formula (6.2), we obtain
$$
1-uz^2+2z^4 g(z^2) = \Bigl(1+uz+2z^2{C(z)\over D(z)}\Bigr) (1-uz+2z^2 g(z)),
$$
or
$ A^*(z)	+ B^*(z) g(z) + C^*(z) g(z^2)=0, $
where
$$
\leqalignno{
	A^*(z)&=(1-uz)C(z)-{u(u-1)\over 2}D(z), \cr
	B^*(z)&=(1+uz)D(z)+2z^2C(z), \cr
	C^*(z)&=-z^2D(z).\cr
}
$$
Since $D(0)=1$, we have $B^*(0)=1$, $C^*(0)=0$, and $C^*(z)\not=0$.
So the power series $g(z)\, ({\rm mod}\, 2)$ satisfies the equation (5.3) with condition~(i).
By Theorem 5.2 (i), the sequence $H\big(g\, ({\rm mod}\, 2)\big)=H(g)\, ({\rm mod}\, 2)$ of Hankel determinants is ultimately periodic
over the field~$\setF_2$.
On the other hand, Identity (6.2) can be rewritten as:
$$
{1\over f(z)}-1+uz = 2 z^2 g(z).
$$
Consequently we obtain
$$
{1\over f(z)}-1+uz \, ({\rm mod}\, 4) = 2 z^2 \times \bigl(g(z)\, ({\rm mod}\, 2)\bigr).\leqno{(6.4)}
$$
Since $f(z)$ is not a rational function modulo $4$ and $f(0)=1$,
the power series $1/f(z)$ is not a rational function modulo $4$. 
Then by Relation (6.4), we know that the power series
$g(z)\,({\rm mod}\, 2)$ is not a rational function. Combining this result with the fact that the sequence
$H(g)\, ({\rm mod}\, 2)$ of Hankel determinants is ultimately periodic
over the field~$\setF_2$, we deduce at once that there exists
an increasing sequence of positive integers
$(n_i)_{i\geq 0}$  such that
$H_{n_i}(g)\not=0$ for all integers $i\geq 0$ and
$\lim\limits_{i\rightarrow \infty} {n_{i+1}\over n_i} =1$.
By Relation~(6.3), we have also $H_{n_i+1}(f)\not=0$ for all integers $i\geq 0$
and $\lim\limits_{i\rightarrow \infty} {n_{i+1}+1\over n_i+1} =1$.
Let $b\geq 2$ be an integer such that for all integers $m\geq 0$, we have
$D({1\over b^{2^m}})f({1\over b^{2^m}})\neq 0$, then $D({1\over b^{2^m}})\neq 0$, and
$$
\Big(1+{{u}\over{b^{2^m}}}\Big)D\Big({1\over b^{2^m}}\Big)+{{2}\over{b^{2^{m+1}}}}C\Big({1\over b^{2^m}}\Big)\neq 0,
$$
hence it follows from Theorem 4.1 that $f(1/b)$ is transcendental and its irrationality exponent
is equal to $2$. \qed

\section{7. Some applications}

\noindent {\it Proof of Theorem 2.3}. Assume $\beta\neq \alpha+1$.
From the definition, we know directly that the power series $F_{\alpha,\beta}(z)$ and $G_{\alpha,\beta}(z)$
converge in the unit disk, and their coefficients in power series expansion are integers.
Moreover we have also
$$
\displaylines{
\rlap{(7.1)}\hfill -1 +  (1+z^{2^\beta})F_{\alpha,\beta}(z)  - z^{2^\alpha}(1+z^{2^\beta})F_{\alpha,\beta}(z^2)= 0.,\hfill\cr
\rlap{(7.2)}\hfill -1 +  (1-z^{2^\beta})G_{\alpha,\beta}(z)  - z^{2^\alpha}(1-z^{2^\beta})G_{\alpha,\beta}(z^2)= 0..\hfill
}
$$
The above equations are of type (6.1). By Theorem 6.1 (i), to conclude,
it suffices to show that $F(z):=F_{\alpha,\beta}(z)\,({\rm mod}\, 2)=G_{\alpha,\beta}(z)\,({\rm mod}\, 2)$ 
is not rational over $\setF_2$. Put
$$
P(t)=z^{2^\alpha}(1+z^{2^\beta})t^2+(1+z^{2^\beta})t+1\in \setF_2(z)[t].
$$
We have $P(F(z))=0$ by (7.1).
By contradiction, suppose that $F(z)$ is rational over $\setF_2$.
Then $P(t)$ is reducible over  $\setF_2(z)$. As a result,
we can find $A(z)$, $B(z)$, and $C(z)$, $D(z)$ in $\setF_2[z]$ such that
$$
P(t)=(A(z)t+B(z))(C(z)t+D(z)).
$$
Then $B(z)D(z)=1$, $A(z)C(z)=z^{2^\alpha}(1+z)^{2^\beta}$, and
$$
A(z)D(z)+B(z)C(z)=1+z^{2^\beta}=(1+z)^{2^\beta},\leqno{(7.3)}
$$
thus $B(z)=D(z)=1$. From the fact that both $z$ and $1+z$ are irreducible over $\setF_2$,
we can find two integers $m,n$ such that $0\leq m\leq 2^\alpha$, $0\leq n\leq 2^\beta$, and
$A(z)=z^{m}(1+z)^{n}$, $C(z)=z^{2^\alpha-m}(1+z)^{2^\beta-n}$. By (7.2), we obtain
$$
z^{m}(1+z)^{n}+z^{2^\alpha-m}(1+z)^{2^\beta-n}=(1+z)^{2^\beta},
$$
from which we deduce necessarily
$$
\leqalignno{
	z^{m}+z^{2^\alpha-m}(1+z)^{2^\beta-2n}&=(1+z)^{2^\beta-n}, \qquad\ \hbox{(if $0\leq n \leq 2^{\beta-1}$),}\cr
	z^{m} (1+z)^{2n-2^\beta}+z^{2^\alpha-m}&=(1+z)^{n}, \qquad \qquad \hbox{(if $2^{\beta-1} < n \leq 2^\beta$).}\cr
}
$$
Put $z=1$ in any one of the above two formulas, we get $2^\beta=2n$,
otherwise the left hand side gives $1$ while the right hand side yields $0$.
Hence, $\beta\geq 1$ and
$z^{m}+z^{2^\alpha-m}=1+z^{2^{\beta-1}}$.
We have either $m=0,\ 2^\alpha-m=2^{\beta-1}$ or
$2^\alpha-m=0,\ m=2^{\beta-1}$.
We deduce at once $\alpha=\beta-1$ in both cases.
This situation has already been excluded,
so the desired result holds. \qed
\vskip 10pt

\noindent {\it Proof of Theorem 2.4}.
It is well known (see [BV13]) that
$$
\leqalignno{
S(z)&=(1+z+z^2)S(z^2),\cr
T(z)&=2-(1+z+z^2)T(z^2). \cr
}
$$
On the other hand, Han has recently shown in [H15b]
that the Hankel determinants of $S(z)$ and $T(z)$ satisfy, for all integers $n\geq 2$,
$$
{H_n(S) \over 2^{n-2}} \equiv {H_n(T) \over 2^{n-2}} \equiv
\cases{
0, \quad \hbox{if $n \equiv 0, 1$ mod $4$}, \cr
1, \quad \hbox{if $n \equiv 2, 3$ mod $4$}. \cr}
$$
Hence there exists
an increasing sequence of positive integers
$(n_i)_{i\geq 0}$  such that
$H_{n_i}(S)\not=0$,
$H_{n_i}(T)\not=0$ for all integers $i\geq 0$ and
$\lim\limits_{i\rightarrow \infty} {n_{i+1}\over n_i} =1$.
It follows from Theorem 4.1 that, for all integers $b\geq 2$, both
$S(1/b)$, $T(1/b)$ are transcendental, and their irrationality
exponents are equal to $2$. \qed

\medskip

We give further concrete examples of transcendental numbers with irrationality exponent equal to $2$.

In [Va15], V\"a\"an\"anen studied the following two power series
$$
L(z)=\sum_{j=0}^\infty {z^{2^j} \over \prod_{i=0}^{j-1} (1-z^{2^i})},
\quad
M(z)=\sum_{j=0}^\infty {(-1)^jz^{2^j} \over \prod_{i=0}^{j-1} (1-z^{2^i})},
$$
which converge in the unit disk with integer coefficients in power series expansion,
and satisfy respectively the functional equations
$$
\leqalignno{
	z(z-1)+(1-z)L(z)-L(z^2)&=0, \cr
	z(z-1)+(1-z)M(z)+M(z^2)&=0. \cr
}
$$
One can check directly that neither $L(z)$ nor $M(z)$ is a rational function modulo 2.
By Theorem 6.1 (ii), we obtain the following result, of which the second part
was proved firstly by V\"a\"an\"anen [Va15].

\proclaim Theorem 7.1.
For all integers $b\geq 2$, both $L(1/b)$ and $M(1/b)$ are transcendental
and their irrationality exponents are equal to $2$.

In a forthcoming paper [FH15], the Hankel determinants of the following power series
$F_5$, $F_{11}$, $F_{13}$, $F_{17a}$ and $F_{17b}$, satisfying the equations
$$
\leqalignno{
	F_5(z)&=(1-z-z^2-z^3+z^4)\, F_5(z^5),\cr
	F_{11}(z)&=(1-z-{z}^{2}+{z}^{3}-{z}^{4}+{z}^{5}+{z}^{6}+{z}^{7}+{z}^{8}-{z}^{9}-{z
}^{10} )\, F_{11}(z^{11}),\cr
F_{13}(z)&=(1-z-{z}^{2}+{z}^{3}-{z}^{4}-{z}^{5}-{z}^{6}-{z}^{7}-{z}^{8}\cr
&\kern 15mm +{z}^{9}-{z}^{10}-{z}^{11}+{z}^{12}) \, F_{13}(z^{13}) \cr
F_{17a}(z)		&=(1-z-z^2+z^3-z^4+z^5+z^6+z^7+z^8+z^9\cr
&\kern 15mm   +z^{10} +z^{11}-z^{12}+z^{13}-z^{14}-z^{15}+z^{16}) \, F_{17a}(z^{17}),\cr
F_{17b}(z)		&=(1-z-z^2-z^3+z^4+z^5-z^6+z^7+z^8+z^9\cr
&\kern 15mm -z^{10} +z^{11}+z^{12}-z^{13}-z^{14}-z^{15}+z^{16})\, F_{17b}(z^{17})\cr
}
$$
are studied and are shown to verify the following relations
$$
\displaylines{
H_n(F_5)/2^{n-1} \equiv H_n(F_{11})/2^{n-1} \equiv  H_n(F_{13})/2^{n-1} \equiv 1 \pmod 2, \cr
H_n(F_{17a})/2^{n-1} \equiv  H_n(F_{17b})/2^{n-1} \equiv 1 \pmod 2.\cr
}
$$
All these power series converge in the unit disk with integer coefficients in power series expansion, and satisfy
the conditions of Theorem 4.1 for all integers $b\geq 2$, thus we obtain

\proclaim Theorem 7.2.
For all integers $b\geq 2$,
all the $F_5(1/b)$, $F_{11}(1/b)$, $F_{13}(1/b)$, $F_{17a}(1/b)$, $F_{17b}(1/b)$
are transcendental and their irrationality exponents are equal to $2$.

\vskip 0.5 cm
\noindent{\bf{Funding.}} This work was supported partially by the National Natural Science Foundation
of China [grant numbers 10990012, 11371210 to\break
J.-Y.Y., 11271223 to Z.-Y.W.]; and the Morningside Center of
Mathematics (CAS) to [J.-Y.Y. and Z.-Y. W.].
\vskip 5pt


\centerline{\bf References}
\bigskip

{
\eightpoint

\article Ad10|Adamczewski, B|On the expansion of some exponential periods in an integer base|Math. Ann.|346|2010|107--116|

\article AB07|Adamczewski, B.;  Bugeaud, Y|On the complexity of algebraic
numbers. I. Expansions in integer bases|Ann. of Math.|165|2007|547--565|

\article AC06|Adamczewski, B.; Cassaigne, J|Diophantine properties of real numbers generated by finite
automata|Compositio Math.|142|2006|1351--1372|

\article AR09|Adamczewski, B.; Rivoal, T|Irrationality measures for some automatic real numbers|
Math. Proc. Cambridge Philos. Soc.|147|2009|659--678|

\article Al87|Allouche, J.-P|Automates finis en th\'eorie des nombres|Exposition. Math.|5|1987|239--266|

\article APWW98|Allouche, J.-P.; Peyri\`ere, J.; Wen, Z.-X.; Wen, Z.-Y|
Hankel determinants of the Thue--Morse sequence|Ann. Inst. Fourier (Grenoble)|48|1998|1--27|

\livre AS03|Allouche, J.-P.; Shallit, J.|Automatic sequences. Theory, applications, generalizations|
Cambridge University Press, Cambridge, {\oldstyle 2003}. xvi+571 pages|

\divers Ba10|Bacher, R|Twisting the Stern sequence, {\tt arxiv.org/abs/1005.5627}, {\oldstyle 2010}, 19 pages|

\livre BG96|Baker, G. A. Jr.; Graves-Morris, P|Pad\'e approximants|
Second edition. Encyclopedia of Mathematics and its Applications, {\bf 59}.
Cambridge University Press, Cambridge, {\oldstyle 1996}. xiv+746 pages|

\article Be94|Becker, P.-G|$k$-regular power series and Mahler-type functional equations|J. Number Theory|49|1994|269-286|

\article BBC15|Bell, J. P.; Bugeaud, Y.; Coons, M|Diophantine approximation of
Mahler numbers|Proc. London Math. Soc.|110|2015|1157--1206|

\livre Br80|Brezinski, C|Pad\'e-type approximation and general orthogonal polynomials|International Series of Numerical Mathematics, {\bf 50}.
Birkh\"au\-ser Verlag, {\oldstyle 1980}. 250 pages|

\article Bu08|Bugeaud, Y|Diophantine approximation and Cantor sets|Math. Ann|341|2008|677--684|

\article Bu11|Bugeaud, Y|On the rational approximation to the Thue--Morse-Mahler numbers|Ann. Inst. Fourier (Grenoble)|61|2011|2065--2076|

\article BKS11|Bugeaud,Y.; Krieger, D.; Shallit, J|Morphic and automatic words: maximal blocks and Diophantine
approximation|Acta Arith.|149|2011|181--199|

\article BV13|Bundschuh, P.; V\"a\"an\"anen, K|Algebraic independence of the generating functions of Stern's
sequence and of its twist|J. Th\'eor. Nombres Bordeaux|25|2013|43--57|

\article Co13|Coons, M|On the rational approximation of the sum of the reciprocals of the Fermat numbers|Ramanujan J.|30|2013|39--65|

\article CV12|Coons, M.; Vrbik, P|An irrationality measure for regular paperfolding numbers|J. Integer Seq.|15|2012|Article 12.1.6|

\article CZ02|Corvaja, P.; Zannier, U|Some new applications of the subspace
theorem|Compositio Math.|131|2002|319--340|

\article De77|Dekking, M|Transcendance du nombre de Thue--Morse|C. R. Acad. Sci. Paris S\'er. A--B|285|1977|A157--A160|

\article Du14|Dubickas, A|On the approximation of the Thue--Morse generating sequence|Bull. 
Math. Soc. Sci. Math. Roumanie|57|2014|59--71|

\article Fa06|Fatou, P|S\'eries trigonom\'etriques et s\'eries de Taylor|Acta Math.|30|1906|335--400|

\divers FH15|Fu, H.; Han, G.-N|Computer assisted proof for Apwenian sequences related to {H}ankel determinants, {\sl preprint}, 30 pages,  {\oldstyle 2015}|

\divers Gr72|Gros, Luc-Agathon-Louis|Th{\'e}orie du baguenodier par un clerc de notaire lyonnais, {\sl Imprimerie d'Aim\'e Vingtrinier}, {\oldstyle 1872}|

\article GWW14|Guo, Y.-J.; Wen, Z.-X.; Wu, W|On the irrationality exponent of the regular paperfolding numbers|Linear Algebra Appl.|446|2014|237--264|

\article H15a|Han, G.-N|Hankel determinant calculus for the Thue--Morse and related sequences|J. Number Theory|147|2015|374-395|

\divers H15b|Han, G.-N|Hankel continued fraction and its applications,
\vskip 0pt
{\tt http://arxiv.org/abs/1406.1593}, {\oldstyle 2014}, 22 pages|    

\livre HKMP13|Hinz, A. M.; Klav{\v{z}}ar, S.; Milutinovi{\'c},
U.; Petr, C|The tower of {H}anoi---myths and maths|Birkh\"auser/Springer Basel AG, Basel, {\oldstyle 2013}|

\article Ke16|Kempner, A. J.|On transcendental numbers|Trans. Amer. Math. Soc.|17|1916|476--482|

\divers Km79|Kmo\v{s}ek, M|Rozwini\c{e}cie niekt\'{o}rych liczb niewymiernych
na u{\l}amki {\l}a\'{n}cu\-chowe. Master's
thesis, Univwersytet Warszawski, {\oldstyle 1979}|

\livre La95|Lang, S|Introduction to Diophantine approximations|Second edition. Sprin\-ger-Verlag, New York, {\oldstyle 1995}. x+130 pages|

\divers LVdP77|Loxton, J. H.; Van der Poorten, A. J|Transcendence and algebraic independence by a method of Mahler.
Transcendence theory: advances and applications (Proc. Conf., Univ. Cambridge, Cambridge, {\oldstyle 1976}), pp. 211--226.
Academic Press, London, {\oldstyle 1977}|

\article Ma29|Mahler, K|Arithmetische Eigenschaften der L\"{o}sungen einer Klass von Funktionalgleichungen|
Math. Ann.|101|1929|342-366|

\article Ma30a|Mahler, K|Arithmetische Eigenschaften einer Klasse transzendental-trans\-zen\-denter
Funktionen|Math. Z.|32|1930|545--585|

\article Ma30b|Mahler, K|\"{U}ber das Verschwinden von Potenzreihen mehrerer Ver\"{a}nderli\-chen in
speziellen Punktfolgen|Math. Ann.|103|1930|573--587|

\article Ni90|Nishioka, K|New approach in Mahler's method|J. Reine Angew. Math.|407|1990|202-219|

\divers Pe14|Pedersen, S. H|A Cantor set type result in the field of formal Laurent series,
{\tt http://arxiv.org/abs/1409.0352}, {\oldstyle 2014}, 12 pages|

\divers Ph15|Philippon, P|Groupes de Galois et nombres automatiques,
\vskip 0pt
{\tt http://arxiv.org/abs/1502.00942}, {\oldstyle 2015}, 22 pages|    

\divers Ro55|Roth, K. F|Rational approximations to algebraic numbers, Mathematika, {\bf 2} ({\oldstyle 1955}), pp. 1--20; corrigendum, 168|

\livre Sa63|Salem, R|Algebraic numbers and Fourier analysis|D. C. Heath and Co., Boston, Mass., {\oldstyle 1963}. x+68 pages|

\article Sh79|Shallit, J|Simple continued fractions for some irrational numbers|
J. Number Theory|11|1979|209--217|

\article SW88|Sharif, H.; Woodcock, C. F|Algebraic functions over a field of positive
characteristic and Hadamard products|J. London Math. Soc.|37|1988|395--403|

\article St58|Stern, M. A|\"Uber eine zahlentheoretische Funktion|J. Reine Angew. Math.|55|1858|193--220|

\divers Va15|V\"a\"an\"anen, K|On rational approximations of certain Mahler functions with a connection to the Thue--Morse sequence,
Int. J. Number Theory, {\oldstyle 2015}, 7 pages, DOI: 10.1142/S1793042115500244|    

\article WW14|Wen, Z.-X; Wu, W|Hankel determinants of the Cantor sequence (Chinese)|
Scientia Sinica Mathematica|44|2014|1059-1072|

\article Za98|Zannier, U|On a functional equation relating a Laurent series
$f(x)$ to $f(x^m)$|Aequations Math.|55|1998|15--43|

\bigskip

}
\goodbreak
\bigskip \noindent
Yann BUGEAUD,
I.R.M.A.,
Universit\'e de Strasbourg et CNRS,
7, rue Ren\'e-Descartes,
67084 Strasbourg, France,
{\tt bugeaud@math.unistra.fr}

\medskip \noindent
Guo-Niu HAN,
I.R.M.A.,
Universit\'e de Strasbourg et CNRS,
7, rue Ren\'e-Descartes,
67084 Strasbourg, France,
{\tt guoniu.han@unistra.fr}

\medskip \noindent
Zhi-Ying WEN,
Department of Mathematics, Tsinghua University, Beijing 100084, China,
{\tt wenzy@tsinghua.edu.cn}

\medskip \noindent
Jia-Yan YAO,
Department of Mathematics, Tsinghua University, Beijing 100084, China,
{\tt jyyao@math.tsinghua.edu.cn}

\end